\documentclass[11pt]{amsart}

\usepackage{latexsym}
\usepackage{amssymb}
\usepackage[cp850]{inputenc}
\usepackage{epsfig}
\usepackage{psfrag}
\usepackage{amsthm}
\usepackage{amscd}
\usepackage{amsmath}
\usepackage{amsfonts}
\usepackage{graphics}
\usepackage[all]{xy}
\usepackage{microtype}

\newtheorem{sat}{Theorem}[section]
\newtheorem{lem}[sat]{Lemma}

\newtheorem{prop}[sat]{Proposition}

\newtheorem*{defi*}{Definition}
\newtheorem*{bei*}{Example}
\newtheorem*{sat*}{Theorem}
\newtheorem*{kor*}{Corollary}
\newtheorem*{rmk*}{Remark}
\newtheorem*{quest*}{Question}
\newtheorem*{fact*}{Fact}

\let\ssection=\section
\renewcommand{\section}{\setcounter{equation}{0}\ssection}

\newtheorem*{namedtheorem}{\theoremname}
\newcommand{\theoremname}{testing}
\newenvironment{named}[1]{\renewcommand{\theoremname}{#1}\begin{namedtheorem}}{\end{namedtheorem}}

\theoremstyle{remark}
\newtheorem*{bem}{Remark}

\newcommand{\BC}{\mathbb C}			\newcommand{\BH}{\mathbb H}
\newcommand{\BR}{\mathbb R}			
			
\newcommand{\BS}{\mathbb S}			\newcommand{\BZ}{\mathbb Z}
			\newcommand{\BT}{\mathbb T}

\newcommand{\CA}{\mathcal A}		
\newcommand{\CC}{\mathcal C}		
		
\newcommand{\CG}{\mathcal G}

		\newcommand{\CT}{\mathcal T}

\newcommand{\actson}{\curvearrowright}
\newcommand{\D}{\partial}

\DeclareMathOperator{\Diff}{Diff}	
\DeclareMathOperator{\PSL}{PSL}		
\DeclareMathOperator{\GL}{GL}		
\DeclareMathOperator{\Id}{Id}		
\DeclareMathOperator{\Map}{Map}
\DeclareMathOperator{\Homeo}{Homeo}
\DeclareMathOperator{\PMap}{PMap}

\newcommand{\abs}[1]{\left\lvert#1\right\rvert}
\newcommand{\z}{\mathbf{z}}

\entrymodifiers={+!!<0pt,\fontdimen22\textfont2>}

\begin{document}

\title[Mapping classes not realized by diffeomorphisms]{Some groups of
  mapping classes\\ not realized by diffeomorphisms}
\author{Mladen Bestvina, Thomas Church \& Juan Souto}
  
\begin{abstract}
  Let $\Sigma$ be a closed surface of genus $g\ge 2$ and $z\in\Sigma$
  a marked point. We prove that the subgroup of the mapping class
  group $\Map(\Sigma,z)$ corresponding to the fundamental group
  $\pi_1(\Sigma,z)$ of the closed surface does not lift to the group
  of diffeomorphisms of $\Sigma$ fixing $z$.
  As a corollary, we show that the Atiyah--Kodaira surface bundles
  admit no invariant flat connection, and obtain another proof of
  Morita's non-lifting theorem.
\end{abstract}
\maketitle

\section{Introduction}
Given a closed orientable surface $\Sigma$ and a finite, possibly
empty, set $\z\subset\Sigma$ of marked points, consider the group
\[\Diff_+(\Sigma,\z)=\{f\in\Diff_+(\Sigma)\vert f(\z)=\z\}\]
of orientation-preserving diffeomorphisms of $\Sigma$ which map the
set of marked points to itself.
(When $\z$ is empty we drop it from our notation.) We denote by
$\Diff_0(\Sigma,\z)$ the normal subgroup of $\Diff_+(\Sigma,\z)$
consisting of those diffeomorphisms which are isotopic to the identity
via an isotopy which fixes the set $\z$. The mapping class group is
the quotient group
\[\Map(\Sigma,\z)=\Diff_+(\Sigma,\z)/\Diff_0(\Sigma,\z).\]
In \cite{Morita}, Morita proved that if $\Sigma$ has genus at least
$18$ and the set of punctures is empty, then the exact sequence
\[0\to\Diff_0(\Sigma) \to\Diff_+(\Sigma)\to\Map(\Sigma)\to 0\] does
not split. The bound was later improved to genus at least 5 by Morita
\cite[Theorem 4.21]{MoritaBook}. Recently Franks--Handel \cite{F-H}
have extended this result so that it holds for genus at least
$3$. Cantat--Cerveau \cite{Cantat} have proved that finite index
subgroups of the mapping class group do not lift to the group of
analytic diffeomorphisms. A much more powerful result is due to
Markovi\'{c} \cite{Markovic} and Markovi\'{c}--\v{S}ari\'{c}
\cite{Markovic-Saric}, who have proved that for genus at least $2$,
the mapping class group does not even lift to the group of
homeomorphisms. The proofs of at least some of these results apply
also to the case with marked points.

Given a subgroup $\Gamma\hookrightarrow \Map(\Sigma,\z)$, the
\emph{realization problem} asks whether $\Gamma$ lifts to
$\Diff_+(\Sigma,\z)$. This has been the focus of much interest for
various classes of subgroups over the years since Nielsen first raised
the question. Affirmative answers were given for cyclic groups by
Nielsen \cite{Nielsen}, for finite groups by Kerckhoff
\cite{Kerckhoff}, and for abelian groups by Birman--Lubotzky--McCarthy
\cite{BLM}. In this paper, we exhibit rather small subgroups of
$\Map(\Sigma,\z)$ that do not lift to $\Diff_+(\Sigma,\z)$.
Specifically, in the case of a surface of genus at least $2$ with a
single marked point we prove:
\begin{sat}\label{weakmeat}
  Let $\Sigma$ be a closed surface of genus $g\ge 2$ and $z\in\Sigma$
  a marked point. No finite index subgroup of the point-pushing
  subgroup ${\pi_1(\Sigma,z)\subset \Map(\Sigma,z)}$ lifts to
  $\Diff_+(\Sigma,z)$.
\end{sat}
The point-pushing subgroup fits into the Birman exact sequence
\begin{equation}\label{birman}
  1\to\pi_1(\Sigma,z)\overset{F}{\to}\Map(\Sigma,z)
  \to\Map(\Sigma)\to 1
\end{equation} as long as $g\geq 2$.
Observe that if $(\Sigma,z)$ is a torus with a single marked point,
then the mapping class group does in fact lift to $\Diff_+(\Sigma,z)$.

We sketch now the proof of Theorem~\ref{weakmeat}. Seeking a
contradiction, assume that there is a homomorphism $\Phi$ such that
the following diagram commutes:
\[\xymatrix{
  & \Diff_+(\Sigma,z) \ar[d]\\
 \pi_1(\Sigma,z) \lhook\mkern-7mu \ar[r]^{F}\ar@{-->}[ru]^{\Phi} &
  \Map(\Sigma,z) }\] 
where $F$ is the inclusion from \eqref{birman}. The homomorphism
$\Phi$ yields an action of $\pi_1(\Sigma,z)$ on $\Sigma$ by
diffeomorphisms fixing $z$ and hence a representation of
$\pi_1(\Sigma,z)$ in $\GL^+(T_z\Sigma)$. By Milnor's inequality this
representation has Euler-number bounded in absolute value by $g-1$. On
the other hand, we compute that the Euler-number must be $2-2g$; this
contradiction gives Theorem~\ref{weakmeat}.

Combining Theorem~\ref{weakmeat} with some topological constructions,
we show that the centralizers of most finite order elements of
$\Map(\Sigma)$ do not lift to $\Diff_+(\Sigma)$. Concretely, we
construct a subgroup of $\Map(\Sigma)$ isomorphic to
$\BZ/3\BZ\times \pi_1(S,z)$ for some closed surface $S$ that does not lift to
$\Diff_+(\Sigma)$. This relies
on the existence of finite order elements and thus does not apply to
finite index subgroups of $\Map(\Sigma)$. Using Theorem~\ref{weakmeat}
and this construction, we derive the following version of Morita's
theorem:

\begin{sat}[Morita's non-lifting theorem]\label{main}
  Let $(\Sigma,\z)$ be a surface of genus $g$ with $\abs{\z}=k$ marked
  points. Assume either that $g\ge 6$ or that $g\ge 2$ and $k\ge
  1$. Then the exact sequence
  \begin{equation}\label{nosplit1}
    0\to\Diff_0(\Sigma,\z)\to\Diff_+(\Sigma,\z)\to\Map(\Sigma,\z)\to 0
  \end{equation}
  does not split. In fact, if $g\ge 2$ and $k\ge 1$ then no finite
  index subgroup of $\Map(\Sigma,\z)$ lifts to $\Diff_+(\Sigma,\z)$.
\end{sat}

Morita originally proved his theorem by finding a surface bundle over
an $6$--dimensional manifold with a cohomological obstruction to the
existence of a flat connection. (All connections are taken to be
smooth.)  The theorem of Earle--Eells \cite{EarleEells} on the
contractibility of $\Diff_0(\Sigma)$ implies that a $\Sigma$--bundle
over a base $B$ admits a flat connection if and only if the
topological monodromy representation $\pi_1(B)\to\Map(\Sigma)$ can be
lifted to a map $\pi_1(B)\to\Diff_+(\Sigma)$. In particular, if the
sequence \eqref{nosplit1} split, then every surface bundle would admit
a flat connection, so Morita's theorem follows from his example.

In contrast, for surface bundles over surfaces, Kotschick--Morita
\cite{KM} proved that every surface bundle admits a flat connection
after ``stabilization''; in particular, there can be no cohomological
obstruction to flatness in this case. This raised the open problem of
finding a surface bundle over a surface that does not admit a flat
connection. The details of the proof of Theorem~\ref{main} give a
partial solution to this problem. In the case of a punctured surface,
Theorem~\ref{weakmeat} gives a surface group isomorphic to
$\pi_1(\Sigma,z)$ inside $\Map(\Sigma,z)$ that does not lift to
$\Diff_+(\Sigma,z)$. This yields a surface bundle with a distinguished
section, with base space a closed surface, which admits no flat
connection such that the distinguished section is parallel. (In fact,
this bundle is just the trivial bundle $\Sigma\times\Sigma$, and the
distinguished section is the diagonal.) We believe that this is the
first such surface group inside a punctured mapping class group
known. In the case of a closed surface, the construction described
above corresponds to a topological construction of Kodaira and Atiyah,
and we conclude (see remarks preceding the proof for definitions):

\begin{sat}\label{AK} When $k\geq 3$,
  the Atiyah--Kodaira bundle $\Sigma\to M_k\to S'$ admits no flat
  connection invariant under the order--$k$ deck transformation
  $\CT\colon M_k\to M_k$.
\end{sat}
However, the full question remains open in the case when the surface
is closed.
\begin{quest*}
  Does there exist a closed surface bundle over a surface that admits
  no flat connection?
\end{quest*}

\noindent {\bf Acknowledgements.} The authors would like to thank
Benson Farb and Vlad Markovi\'{c} for their interest in this
project. The second author would like to thank Benson Farb for
introducing him to the examples of Kodaira and Atiyah and to the
questions surrounding flat surface bundles. We are very grateful to an
anonymous referee for their careful reading, and for pointing out that
the bound in our main theorem could be improved.

\section{A few facts about Euler-numbers}
Let $\Sigma$ be a closed surface of genus $g$ and let
$\widetilde\Sigma\to\Sigma$ be its universal cover. Choose base points
$z\in\Sigma$ and $\tilde z\in\widetilde\Sigma$ projecting to $z$. The
choice of base points yields an identification between the fundamental
group $\pi_1(\Sigma,z)$ and the deck-transformation group of the
cover $\widetilde\Sigma\to\Sigma$. Before going any further, let us
remark that the composition $\gamma\star\eta$ of two elements
$\gamma,\eta\in\pi_1(\Sigma,z)$ is obtained by first running $\gamma$
and then $\eta$. By construction, the universal cover
$\widetilde\Sigma$ consists of homotopy classes rel endpoints of
continuous paths in $\Sigma$ beginning at $z$. Here we can identify
$\tilde z$ with, for instance, the homotopy class of the constant
path. The fundamental group $\pi_1(\Sigma,z)$ acts on
$\widetilde\Sigma$ by precomposition, meaning that we first run a path
representing the element in the fundamental group and then a path
representing the element in $\widetilde\Sigma$. In particular, the
obtained action of $\pi_1(\Sigma,z)\actson\widetilde\Sigma$, the
so-called action by deck-transformations, is a left action.

Assume now that $\rho\colon\pi_1(\Sigma,z)\to\Homeo^+(\BS^1)$ is an
action of the fundamental group of $\Sigma$ on the circle. Let
$E_\rho$ be the quotient of $\widetilde\Sigma\times\BS^1$ under the action
\[\pi_1(\Sigma,z)\actson(\widetilde\Sigma\times\BS^1),\quad
(\gamma,(x,\theta))\mapsto(\gamma x,\rho(\gamma)\theta).\]

\noindent The projection of $\widetilde\Sigma\times\BS^1$ onto the
first factor is $\pi_1(\Sigma)$--equivariant and has fiber $\BS^1$; this descends to
give $E_\rho$ the structure of a circle bundle over $\Sigma$. The
trivial connection on $\widetilde\Sigma\times \BS^1$ induces a flat
connection on $E_\rho$. Conversely, every flat circle bundle over
$\Sigma$ is obtained in this way.

The \emph{Euler-number} $e(E_\rho)\in\BZ$ of the bundle $E_\rho\to\Sigma$ is
the obstruction for the bundle $E_\rho$ to admit a section, or
equivalently, for the action $\rho$ to lift to an action on the
universal cover $\BR$ of $\BS^1$.

\begin{named}{Milnor--Wood inequality}
  Assume that $E_\rho$ is a flat orientable circle bundle over a
  closed surface $\Sigma$ of genus $g$. Then $\abs{e(E_\rho)}\le
  2g-2$.
\end{named}

It should be observed that there are flat circle bundles with
Euler-number $2-2g$. For instance, endowing $\Sigma$ with a hyperbolic
metric, we can identify the universal cover $\widetilde\Sigma$ with the
hyperbolic plane. The action of $\pi_1(\Sigma,z)$ on $\BH^2$ extends
to an action on the circle at infinity $\D_\infty\BH^2$. The
associated flat circle bundle is isomorphic to the unit tangent bundle
of $\Sigma$ and hence has Euler-number equal to the
Euler characteristic $\chi(\Sigma)=2-2g$. We record this
fact for further reference (see \cite[Appendix C]{MS}):

\begin{lem}\label{unit-tangent}
  Let $\Sigma$ be a closed orientable hyperbolic surface of genus $g$
  and identify $\pi_1(\Sigma,z)$ with the corresponding group of
  deck-transformations of $\BH^2$. The circle bundle corresponding to
  the induced action of $\pi_1(\Sigma,z)$ on $\D_\infty\BH^2=\BS^1$
  has Euler-number $2-2g$.\qed
\end{lem}

We point out that Goldman \cite{G} proved a converse to this lemma: if
$\rho\colon\pi_1(\Sigma,z)\to \PSL_2\BR$ has $\abs{e(E_\rho)}=2g-2$,
then $\rho$ is an isomorphism onto a discrete subgroup of $\PSL_2\BR$
and thus comes from a hyperbolic metric on $\Sigma$ as in the lemma.

Other examples of circle bundles over $\Sigma$ can be constructed as
follows. A linear action $\rho\colon\pi_1(\Sigma,z)\to\GL_2^+\BR$ of
$\pi_1(\Sigma,z)$ on $\BR^2$ induces an action on the space of
directions $P_+\BR^2=(\BR^2\setminus\{0\})/\BR_+$ of $\BR^2$. The latter can be
identified with the circle and hence the same construction as above
yields a circle bundle $E_\rho$.  A circle bundle $E_\rho$ arising in
this way is called a \emph{flat linear circle bundle}. The linear
action $\rho$ induces a different circle bundle $\hat E_\rho$ via the
induced projective action on the projective line
$P\BR^2=(\BR^2\setminus\{0\})/(\BR\setminus\{0\})$, which can also be
identified with the circle.  By construction there is a two-to-one
fiberwise covering $E_\rho\to\hat E_\rho$. In particular, $e(\hat
E_\rho)=2e(E_\rho)$. We have then:

\begin{named}{Milnor's inequality}
  Assume that $E_\rho$ is a flat linear orientable circle bundle over
  a closed surface $\Sigma$ of genus $g$. Then $\abs{e(E_\rho)}\le
  g-1$.
\end{named}

In \cite{Milnor}, Milnor proved that if a $\GL_2^+\BR$--bundle over a
closed surface of genus $g$ admits a flat symmetric connection, then
its Euler-number is bounded in absolute value by $g-1$. This is
equivalent to Milnor's inequality above. Later, Wood \cite{Wood}
extended Milnor's work to prove the Milnor--Wood inequality.

For a general oriented circle bundle $S^1\to E \to B$, the Euler class is a
characteristic class $e(E)\in H^2(B)$. When the base space is a
surface, we identify this with the Euler-number by the identification
$H^2(\Sigma)= \BZ$. We will use the same symbol for the Euler-number
and Euler class; it should be clear from context what is meant.

\section{Surfaces with one puncture}
Let $\Sigma$ be a closed surface of genus $g$ and $z\in\Sigma$ a
marked point, and define the group $\CG(\Sigma,z)$ to consist of
those orientation-preserving homeomorphisms $f$ of $\Sigma$ which fix
$z$ so that $f$ and $f^{-1}$ are differentiable at $z$. In this section we
prove the following generalization of Theorem~\ref{weakmeat}:

\begin{prop}\label{meat}
  Let $\Sigma$ be a closed surface of genus $g\ge 2$ and $z\in\Sigma$
  a marked point. If $\Gamma\subset\pi_1(\Sigma,z)$ is a finite index
  subgroup, then the inclusion of $\Gamma$ into $\Map(\Sigma,z)$ under
  the homomorphism $F$ from \eqref{birman} does not lift to
  $\CG(\Sigma,z)$.
\end{prop}

Observe that since $\Diff_+(\Sigma,z)$ is a subgroup of
$\CG(\Sigma,z)$, Theorem~\ref{weakmeat} follows directly from
Proposition~\ref{meat}. Although Proposition~\ref{meat} applies only
to punctured surfaces, we will upgrade it in Section~\ref{generalcase}
to prove Theorem~\ref{main} for closed surfaces.

Before going any further we describe the homomorphism
\[F\colon\pi_1(\Sigma,z)\hookrightarrow\Map(\Sigma,z)\] from \eqref{birman} in
detail. Given $\gamma\in\pi_1(\Sigma,z)$, let
$\vec\gamma\colon[0,1]\to\Sigma$ be a loop in the corresponding homotopy
class. The map $t\mapsto\vec\gamma(1-t)$ can be interpreted as an
isotopy from the identity $\Id_z$ to itself. By the theorem on
extension of isotopies we obtain an isotopy $f_t\colon\Sigma\to\Sigma$ with
$f_0=\Id_\Sigma$ and $f_t(z)=\vec\gamma(1-t)$. Birman proved that the element
$F_{\gamma}\in\Map(\Sigma,z)$ corresponding to
$f_1\in\Diff_+(\Sigma,z)$ depends only on the element
$\gamma\in\pi_1(\Sigma,z)$. Observing that
\[F_{\gamma\star\eta}=F_{\gamma}\circ F_{\eta}\] we have that
$F\colon\pi_1(\Sigma,z)\to\Map(\Sigma,z)$ is a homomorphism.

Starting now the proof of Proposition~\ref{meat}, assume that there
is a homomorphism
\[\Phi\colon\pi_1(\Sigma,z)\to\CG(\Sigma,z)\]
such that for each $\gamma\in\pi_1(\Sigma,z)$ the homeomorphism
$\Phi_\gamma$ represents the mapping class
$F_\gamma\in\Map(\Sigma,z)$. Endowing $\Sigma$ with a hyperbolic
metric we identify its universal cover with $\BH^2$; choose a point
$\tilde z$ covering $z$. We obtain then a homomorphism
\[\tilde\Phi\colon\pi_1(\Sigma,z)\to\CG(\BH^2,\tilde z)\]
mapping $\gamma$ to the unique lift of $\Phi_\gamma$ which fixes
$\tilde z$. Here $\CG(\BH^2,\tilde z)$ is the group of homeomorphisms
of $\BH^2$ fixing $\tilde z$ which are differentiable at $\tilde z$
with inverse differentiable at $\tilde z$.

\begin{lem}\label{extension-fix}
  The homeomorphism $\tilde\Phi_\gamma\colon\BH^2\to\BH^2$ extends to a
  homeomorphism of the closed disk
  $\overline{\BH}^2=\BH^2\cup\D_\infty\BH^2$. Moreover, the
  restriction of $\tilde\Phi_\gamma$ to $\D_\infty\BH^2$ coincides
  with the action of $\gamma$ as a deck-transformation.
\end{lem}

Lemma~\ref{extension-fix} is probably well-known to experts and
non-experts alike. However, here is a proof:

\begin{proof}
  We start by observing that the action $\Phi$ can be lifted in a
  different way. By construction, if we forget the marked point, the
  homeomorphism $\Phi_\gamma$ is homotopic to the identity. If $f_t$
  is such a homotopy with $f_0=\Id_\Sigma$ and $f_1=\Phi_\gamma$, let
  $\hat{f}_t$ be the unique lift of $f_t$ to $\BH^2$ with
  $\hat{f}_0=\Id_{\BH^2}$. We obtain a new lift
  $\hat\Phi_\gamma=\hat{f}_1$ of $\Phi_\gamma$. It follows directly
  from the construction of the homomorphism $F$ and from the fact that
  $\Phi_\gamma$ represents $F(\gamma)$ that
  \[\hat\Phi_\gamma(\tilde z)=\gamma^{-1}\tilde z\]
  where we have identified $\gamma\in\pi_1(\Sigma,z)$ with the
  corresponding deck-transformation. In particular, the two lifts
  $\hat\Phi_\gamma$ and $\tilde\Phi_\gamma$ differ by the
  deck-transformation $\gamma$, meaning that
  \begin{equation}\label{relation-lifts}
    \gamma\circ\hat\Phi_\gamma=\tilde\Phi_\gamma.
  \end{equation}
  By construction, the lift $\hat\Phi_\gamma$ moves every point in
  $\BH^2$ a uniformly bounded distance from itself. In particular
  $\hat\Phi$ extends continuously to the identity map on the boundary
  $\D_\infty\BH^2$ of the hyperbolic plane. The claim follows from
  this fact and \eqref{relation-lifts}.
\end{proof}

We come now to the meat of the proof of Theorem~\ref{main}. Recall
that $\overline{\BH}^2$ is the union of $\BH^2$ with the circle at
infinity $\partial_\infty \BH^2$. The half-open annulus
$\overline{\BH}^2\setminus\tilde z$ can be compactified in a canonical
way by attaching to the open end the space of directions $P_+T_{\tilde
  z}\BH^2=(T_{\tilde z}\BH^2\setminus\{0\})/\BR_+$ of the tangent
space at $\tilde z$. Let $\CA$ be the so-obtained closed annulus. By
Lemma~\ref{extension-fix}, the action of $\pi_1(\Sigma,z)$ via
$\tilde\Phi$ induces an action on $\overline{\BH}^2\setminus\{\tilde
z\}$. Moreover, the assumption that $\tilde\Phi_\gamma$ is
differentiable at $\tilde z$ for all $\gamma\in\pi_1(\Sigma,z)$
implies that this action extends to an action on $\CA$ which restricts
to $\D\CA$ as follows.
\begin{itemize}
\item On the component $\D_1\CA$ corresponding to $\D_\infty\BH^2$ the
  action of $\pi_1(\Sigma,z)$ the action is equal to the one induced
  by the deck-transformation group by Lemma~\ref{extension-fix}.
\item On the component $\D_2\CA$ corresponding to the space of
  directions of $T_{\tilde z}\BH^2$, the action is induced by the
  representation
  \[\pi_1(\Sigma,z)\to\GL^+(T_{\tilde z}\BH^2),\ \ \gamma\mapsto
  d\hat\Phi_\gamma\vert_{\tilde z} 
  \]
\end{itemize}
In particular, it follows from Lemma~\ref{unit-tangent} that the
circle bundle $E_1$ over $\Sigma$ induced by the action on $\D_1\CA$
has Euler-number
\[e(E_1)=2-2g.\] Similarly, it follows from Milnor's inequality that
the circle bundle $E_2$ over $\Sigma$ induced by the action on
$\D_2\CA$ satisfies
\[\left\vert e(E_2)\right\vert=g-1.\]
But since the annulus bundle $\CA$ admits a fiberwise deformation
retract onto $E_1$ and also onto $E_2$, these bundles have the
same Euler-number
\[e(E_1)=e(\CA)=e(E_2).\] This contradiction shows that the image of
$\pi_1(\Sigma,z)$ under $F$ does not lift to $\CG(\Sigma,z)$. The same
argument applies to finite index subgroups; this concludes the proof
of Proposition~\ref{meat}.\qed
\medskip

As mentioned above, Theorem~\ref{weakmeat} follows directly from
Proposition~\ref{meat}.

\subsection*{An alternate perspective on Proposition~\ref{meat}}
In the remainder of this section, we sketch an alternate perspective
on the above proof in the language of surface bundles. This
perspective will be used in the remarks following the proof of
Theorem~\ref{main} and in the proof of Theorems~\ref{AK} and \ref{m-const}.

The previous section considered the flat linear circle bundle
$E_{d\Phi}\to \Sigma$, which \textit{a priori} depends on the lift
$\Phi$ of $F$; however, the isomorphism type of $E_{d\Phi}$ as a
topological circle bundle does not depend on $\Phi$. In fact, this
circle bundle can be defined without reference to any lift, as we
describe below.

The theorem of Earle--Eells, extended to punctured surfaces by
Earle--Schatz \cite{EarleSchatz}, gives a one-to-one correspondence
between $\Sigma$--bundles with distinguished section over a base $B$
(up to isomorphism) and their monodromy representation
$\pi_1(B)\to\Map(\Sigma,z)$ (up to conjugacy). The ``vertical Euler
class'' of a $\Sigma$--bundle with distingushed section is a
characteristic class defined as follows. Given such a bundle
$\Sigma\to E\overset{\pi}{\to} B$ with section $\sigma\colon B\to E$,
the vectors tangent to the fibers span a 2--dimensional subbundle
$T\pi\leq TE$. Passing to the space of directions and restricting to
the section $\sigma$ induces a circle bundle $UT\pi|_\sigma\to B$. The
vertical Euler class is defined to be the Euler class
$e(UT\pi|_\sigma) \in H^2(B)$ of this circle bundle. This class is
discussed in many references, including \cite{Morita}. We will need
only the following property.

\begin{fact*} If the monodromy $r\colon\pi_1(B)\to\Map(\Sigma,z)$ of a
  $\Sigma$--bundle with section lifts to $\rho\colon
  \pi_1(B)\to\CG(\Sigma,z)$, yielding as above the flat linear circle
  bundle ${E_{d\rho}\to B}$, then $E_{d\rho}$ is isomorphic to
  $UT\pi|_\sigma$ as a circle bundle.
\end{fact*}

To apply this fact to the map $F\colon \pi_1(\Sigma,z)\to
\Map(\Sigma,z)$, we must identify the $\Sigma$--bundle with section
over $\Sigma$ whose monodromy is $F$. It is easy to check that the
desired bundle is the product bundle $p_1\colon \Sigma\times\Sigma\to
\Sigma$, with section given by the diagonal $\Delta\colon \Sigma\to
\Sigma\times \Sigma$.

Along the diagonal, we can identify the tangent space
$T_{(p,p)}(\Sigma\times \Sigma)$ with $T_p\Sigma\times
T_p\Sigma$. Under this identification, $Tp_1=\ker dp_1$ consists of
vectors of the form $(0,v)\in T_p\Sigma\times T_p\Sigma$. Mapping
$(0,v)\mapsto (v,v)$ gives an isomorphism between $Tp_1|_\Delta$ and
$T\Delta$, the subbundle spanned by vectors tangent to the
diagonal. It follows that $e(UTp_1|_\Delta)=e(UT\Delta)=2-2g$. By
Milnor's inequality, this bundle is not isomorphic to any flat linear
circle bundle. Thus the fact above implies that no lift $\Phi\colon
\pi_1(\Sigma,z)\to\CG(\Sigma,z)$ exists.

For a finite index subgroup of $\pi_1(\Sigma,z)$ corresponding
to the cover ${p\colon\Sigma'\to \Sigma}$, the same argument applies to
the bundle $\Sigma'\times \Sigma\to \Sigma$, with section given by the
graph of $p$.

\section{The proof of Theorem~\ref{main}}
\label{generalcase}
In this section we deduce Theorem~\ref{main} from
Proposition~\ref{meat}, but before doing so we need some notation.

\begin{named}{Theorem \ref{main}}
  Let $(\Sigma,\z)$ be a surface of genus $g$ with $k$ marked
  points. Assume that either $g\ge 6$ or that $g\ge 2$ and $k\ge
  1$. Then the exact sequence
  \[0\to\Diff_0(\Sigma,\z)\to\Diff_+(\Sigma,\z)\to\Map(\Sigma,\z)\to
  0\] does not split. In fact, if $g\ge 2$ and $k\ge 1$ then no finite
  index subgroup of $\Map(\Sigma,\z)$ lifts to $\Diff_+(\Sigma,\z)$.
\end{named}

Given a surface as in Theorem~\ref{main}, let $\CG(\Sigma,\z)$ be the
group of those orientation-preserving homeomorphisms $f$ of $\Sigma$
which fix the marked points $\z$ pointwise so that $f$ and $f^{-1}$
are differentiable at each $z\in\z$. If $\CG_0(\Sigma,\z)$ denotes the
normal subgroup of $\CG(\Sigma,\z)$ consisting of those elements which
are isotopic to the identity relative to the set $\z$ then the
quotient group
\[\PMap(\Sigma,\z)=\CG(\Sigma,\z)/\CG_0(\Sigma,\z)\]
is the \emph{pure mapping class group}, a finite index subgroup of the
mapping class group $\Map(\Sigma,\z)$. We could equivalently define
$\PMap(\Sigma,\z)$ using diffeomorphisms instead of $\CG(\Sigma,\z)$.

We can now start with the proof of Theorem~\ref{main}. We will divide
the proof into cases depending on the genus $g$ and number of marked
points $k$ in $(\Sigma,\z)$; the proof for each case will depend upon
the previous one.  \medskip

\noindent \textbf{Case 1.} $g\ge 2$ and $k=1$.
\medskip

Since the group $\Diff_+(\Sigma,z)$ is a subgroup of $\CG(\Sigma,z)$,
the claim follows directly from Proposition~\ref{meat}.\qed\medskip

\noindent \textbf{Case 2.} $g\ge 2$ and $k\ge 2$.
\medskip

Consider the configuration space
\[\CC_k(\Sigma)=
\big\{(x_1,\dots,x_k)\in\Sigma^k\big\vert x_i\neq x_j\ \hbox{if}\ i\neq j\big\}\]
of ordered $k$--tuples of pairwise distinct points in the closed surface
$\Sigma$. We can consider $\CC_k(\Sigma)$ as a fiber bundle over $\Sigma$
via the following projection:
\[p_1\colon\CC_k(\Sigma)\to\Sigma,\ \ p_1\colon(x_1,\dots,x_k)\mapsto x_1\]
In particular, we obtain a homomorphism
\[\pi_1(p_1)\colon\pi_1(\CC_k(\Sigma),(z_1,\dots,z_k))\to\pi_1(\Sigma,z_1).\]
We claim that $\pi_1(p_1)$ has a right inverse:
\begin{lem}\label{tom}
There is a homomorphism 
\[\eta\colon\pi_1(\Sigma,z_1)\to\pi_1(\CC_k(\Sigma),(z_1,\dots,z_k))\]
with $\pi_1(p_1)\circ\eta=\Id$.
\end{lem}
\begin{proof}
  It suffices to construct a section $\Sigma\to\CC_k(\Sigma)$ of the
  fiber bundle $p_1\colon\CC_k(\Sigma)\to\Sigma$. In order to
  construct such a section, it suffices to find maps $\alpha_i\colon
  \Sigma\to \Sigma$ for $i=2,\dotsc,k$, each without fixed points and
  satisfying $\alpha_i(z_1)=z_i$ and $\alpha_i(x)\neq \alpha_j(x)$ for
  $i\neq j$. Given such $\alpha_i$, let
  $\sigma\colon\Sigma\to\Sigma^k$ be the map given by
  $\sigma(x)=(x,\alpha_2(x),\dotsc,\alpha_k(x))$. By construction, the
  image of $\sigma$ is contained in $\CC_k(\Sigma)$. On the other
  hand, $p_1\circ\sigma=\Id$; in other words, $\sigma$ is the desired
  section.

  To find such maps, let $T\subset\Sigma$ be a compact subsurface
  homeomorphic to a torus with one boundary component and which
  contains all the points $z_1,\dots,z_k$. Let $C$ be a homotopically
  essential simple closed curve in $T\setminus\D T$ with $z_i\in C$
  for $i=1,\dotsc,k$; let also $\BT$ be the closed torus obtained by
  collapsing the boundary of $T$ to a point. Equivalently, $\BT$ is
  obtained by collapsing $\Sigma\setminus(T\setminus\D T)$ to a point;
  this gives a map $\Sigma\to\BT$. We can now identify $C$ with a
  factor of $\BT\approx \BS^1\times \BS^1$, giving in particular a
  projection $\BT\to C$. Composing with the map $\Sigma\to\BT$ above,
  we obtain a retraction $a\colon\Sigma\to C$ which fixes each point
  in $C$. Fixing a parametrization of $C$, let $\alpha_i$ be the
  composition
  \[\alpha_i\colon\quad\Sigma\overset{a}{\twoheadrightarrow}C
  \overset{r_i}{\longrightarrow} C\hookrightarrow \Sigma\] where the
  middle map $r_i\colon C\to C$ is the rotation taking $z_1$ to
  $z_i$. Since the image of each $\alpha_i$ is $C$, any fixed point of
  $\alpha_i$ must lie in $C$; since $\alpha_i$ acts by a nontrivial
  rotation on $C$, $\alpha_i$ has no fixed points. Similarly, since
  each $\alpha_i$ is the composition of $a$ with a different rotation,
  we have $\alpha_i(x)\neq \alpha_j(x)$ for $i\neq j$, as desired.
\end{proof}

Order now the points $z_1,\dots,z_k$ in $\z$ and let $\vec\z$ be the
so-obtained point in $\CC_k(\Sigma)$. Recall that $\PMap(\Sigma,\z)$
is the pure mapping class group of $(\Sigma,\z)$, i.e. the subgroup of
the mapping class group consisting of mapping classes whose
representatives in $\Diff_+(\Sigma)$ fix each one of the marked
points. Forgetting all the marked points, and forgetting all the
marked points but $z_1$, we obtain the following versions of the Birman
exact sequence \eqref{birman}:
\[\xymatrix{
1 \ar[r] & \pi_1(\CC^k(\Sigma),\vec\z)\ar[r]\ar[d]^{\pi_1(p)}
        & \PMap(\Sigma,\z)\ar[r]\ar[d] &\Map(\Sigma)\ar@{=}[d]\ar[r] & 1 \\
1 \ar[r] & \pi_1(\Sigma,z_1)\ar@/^/[u]^\eta\ar[r]
        & \Map(\Sigma,z_1)\ar[r] &\Map(\Sigma) \ar[r]  & 1}\]
Here $\eta$ is the homomorphism provided by Lemma~\ref{tom}.

Assume now that $G$ is a finite index subgroup in $\Map(\Sigma,\z)$
which lifts to $\Diff_+(\Sigma,\z)$. Intersecting with the
point-pushing subgroup $\pi_1(\CC^k(\Sigma),\vec\z)$, we obtain a
finite index subgroup of $\pi_1(\CC^k(\Sigma,\vec\z))$ which lifts to
$\Diff_+(\Sigma,\z)$. Composing with the section $\eta$ provided by
Lemma~\ref{tom}, we obtain a lift of a finite index subgroup
$\Gamma<\pi_1(\Sigma,z_1)$ to $\Diff_+(\Sigma,\z)$. Since
$\Diff_+(\Sigma,\z)$ is a subgroup of $\Diff_+(\Sigma,z_1)$ and hence
of $\CG(\Sigma,z_1)$, this contradicts Proposition~\ref{meat}. This
concludes the proof of Case~2.\qed

\begin{bem}
  Before going further, observe that we have actually proved that,
  under the assumptions of Case~2, no finite index subgroup of
  $\Map(\Sigma,\z)$ lifts to $\CG(\Sigma,\z)$.
\end{bem}\bigskip

\textbf{Case 3.} $g\ge 6$ and $k=0$.
\medskip

In this case we will prove that the centralizer of a certain finite
order element $T\in\Map(\Sigma)$ does not lift to
$\Diff_+(\Sigma)$. We have significant freedom in our choice of $T$;
we require only that the order of $T$ be at least 3, and that the
quotient $\Sigma/\langle T\rangle$ have genus at least 2. The first
step is to verify that such finite order elements exist for all
$\Sigma$.  Though in the proof we work with an order 3 automorphism
$\tau$, any number $k\geq 3$ would work just as well;
see the remark following the proof of Lemma~\ref{push} to see why it
is necessary that $\tau$ have order at least 3. \medskip

\begin{fact*} If $g\geq 6$, then there is a diffeomorphism $\tau\colon
  \Sigma\to\Sigma$ of order $3$ with at least 2 fixed points so that
  the quotient $\Sigma/\langle\tau\rangle$ has genus $h\geq 2$.
\end{fact*}
There are many different ways to find such a finite-order
diffeomorphism. One uniform way is to begin with a degree $3$ cyclic
branched cover of the sphere branched at $g-4$ points. By the Hurwitz
formula, the resulting surface has genus $g-6$. Now add three genus 2
handles symmetrically, so they are permuted freely by the order 3 deck
transformation; in the quotient this corresponds to adding a single
genus 2 handle to the sphere. The result is a genus $g$ surface
$\Sigma$ with an order 3 automorphism $\tau$ so that the quotient
$\Sigma/\langle\tau\rangle$ has genus 2.

Let $\tau\colon\Sigma\to\Sigma$ be the diffeomorphism provided by the
fact above, ${T\in\Map(\Sigma)}$ the corresponding mapping
class, and
\[C(T)=\{f\in\Map(\Sigma)\vert
f\circ T=T\circ f\}\] its centralizer. We claim that
$C(T)$ does not lift to $\Diff_+(\Sigma)$. Seeking a
contradiction, assume that such a lifting
\[\Psi\colon C(T)\to\Diff_+(\Sigma)\]
exists. By definition, the diffeomorphism $\Psi(T)$ has order $3$ and
is isotopic to $\tau$. In particular, both diffeomorphisms are
conjugate and we may assume without loss of generality that
$\Psi(T)=\tau$, so that the image of $\Psi$ is contained in the
centralizer $C(\tau)<\Diff_+(\Sigma)$.

\begin{bem}
  The authors did not find a reference for this fact, so we give a
  short argument here. Each of $\tau$ and $\tau'=\Psi(T)$ is an
  isometry of some hyperbolic structure $X$ and $X'$ on $\Sigma$,
  respectively. Identifying the universal cover of $X$ and $X'$ with
  the hyperbolic plane, we obtain that the groups $G$ generated by all
  lifts of $\tau$ and $G'$ generated by all lifts of $\tau'$ are
  Fuchsian groups. In fact, the assumption that $\tau$ is isotopic to
  $\tau'$ implies that $G$ and $G'$ are isomorphic. Satz IV.10 in
  Zieschang--Vogt--Coldewey \cite{ZVC} implies that the actions of $G$
  and $G'$ are conjugate. This yields a conjugation between $\tau$ and
  $\tau'$. Before moving on, we observe that a second and slightly
  more sophisticated proof follows from the fact that the fixed point
  set of the mapping class $T$ in Teichm\"uller space is totally
  geodesic with respect to the Teichm\"uller metric, and thus \textit{a
    fortiori} connected.
\end{bem}

By construction, the quotient surface $S=\Sigma/\langle\tau\rangle$
has genus $h\ge 2$. Let $z_1,\ldots,z_k\in S$ be the projection to $S$
of the fixed points of $\tau$ and set $\z=\{z_1,\ldots,z_k\}$. Every
$f\in\Diff_+(\Sigma)$ which commutes with $\tau$ induces a
homeomorphism of $(S,\z)$. This gives a homomorphism
\[\alpha\colon C(\tau)\to\Homeo(S,\z)\] whose kernel is the cyclic
group generated by $\tau$. Let $C(\tau,\z)$ be the finite index
subgroup of $C(\tau)$ consisting of those diffeomorphisms which commute with
$\tau$ and fix each of its fixed points. The key fact, and the reason
we require $\tau$ to have order 3, is the following lemma:
\begin{lem}\label{push}
  The image of $C(\tau,\z)$ under $\alpha$ is contained in $\CG(S,\z)$.
\end{lem}
\begin{proof}
  It is well-known that there is a conformal structure on $\Sigma$
  such that $\tau$ is biholomorphic. In particular, if $x$ is one of
  the fixed points of $\tau$ we can find coordinates $\zeta$ around
  $x$ such that $\tau(\zeta)=\omega\cdot\zeta$ where $\omega$ is a
  primitive third root of unity. Every differentiable
  $f\colon\Sigma\to\Sigma$ which fixes $x$ and commutes with $\tau$
  has differential
  \[df_x\colon T_x\Sigma\to T_x\Sigma\] satisfying
  $df_x\cdot\omega=\omega\cdot df_x$. Since $\omega$ has order $3$,
  the elements $1$ and $\omega$ span $\BC$ as a real vector
  space. Since $df_x$ commutes with multiplication by each, $df_x$ is
  complex differentiable. This implies that the induced map $S\to S$
  is also differentiable at the projection of $x$. This concludes the
  proof of the lemma. Note that we could not have concluded that
  $df_x$ is complex differentiable if $\omega$ instead had order 2,
  since any linear map commutes with $-1$.
\end{proof}

By composing with $\Psi$, we obtain an action
\[C(T)\overset{\Psi}{\longrightarrow}C(\tau)
\overset{\alpha}{\longrightarrow}\Homeo(S,\z)\] of $C(T)$ on
$(S,\z)$. Since $\langle \tau\rangle$ is the kernel of $\alpha$, this
descends to an action
\[C(T)/\langle T\rangle\to \Homeo(S,\z).\]

As in the construction of $\alpha$, we can identify $C(T)/\langle
T\rangle$ with a certain subgroup of $\Map(S,\z)$. A mapping class in
$\Map(S,\z)$ lifts to the branched cover $\Sigma$ exactly if it
preserves up to conjugacy the subgroup of $\pi_1(S\setminus \z)$
determining the cover $\Sigma\setminus\z\to S\setminus\z$. Since this
subgroup has finite index in $\pi_1(S\setminus \z)$, its stabilizer
has finite index in $\Map(S,\z)$. Among these, $C(T)/\langle T\rangle$
is identified with the finite index subgroup consisting of those
mapping classes whose lift to $\Sigma$ commutes with $T$. Let $\Gamma$
be the intersection of $C(T)/\langle T\rangle$ with $\PMap(S,\z)$;
note that $\Gamma$ has finite index in $\Map(S,\z)$.

We consider the restriction of the action $C(T)/\langle T\rangle\to
\Homeo(S,\z)$ above to the subgroup $\Gamma$. Since $\Gamma$ is
contained in $\PMap(S,\z)$, the image under $\Psi$ of any lift will be
contained in $C(\tau,\z)$. Lemma~\ref{push} implies that the action
$\Gamma\to \Homeo(S,\z)$ has image contained in $\CG(S,\z)$. Thus we
have a lift of the finite index subgroup $\Gamma<\Map(S,\z)$ to
$\CG(S,\z)$, contradicting the remark following the proof of
Case~2. This contradiction completes the proof of Case~3, and thus
concludes the proof of Theorem~\ref{main}.  \qed \bigskip

\noindent For a minimal example of a non-lifting subgroup, consider
the intersection of $\Gamma\subset \Map(S,\z)$ with the surface group
$\eta(\pi_1(S,z_1))$; this gives a surface group inside $\Map(S,\z)$
whose preimage in $C(T)$ does not lift to $\Diff_+(\Sigma)$. This
preimage is a central extension of a surface group by the cyclic group
$\langle T\rangle$; by possibly passing to an index 3 subgroup, we
may assume this extension is trivial, yielding a subgroup of
$\Map(\Sigma)$ isomorphic to $\BZ/3\BZ\times \pi_1(S',z)$ which does
not lift to $\Diff_+(\Sigma)$.

\subsection*{Observations on the proof of Theorem~\ref{main}}
In this section, we give an informal discussion interpreting the above
proof in terms of surface bundles. We then use this perspective to
give two observations, Theorems~\ref{AK} and \ref{m-const} below.

As discussed in the introduction, Case~1 above is equivalent to the
statement that not every surface bundle with section admits a flat
connection so that the section is parallel. This was proved in
Proposition~\ref{meat} by exhibiting the product bundle
$\Sigma\times\Sigma$ with section given by the diagonal $\Delta$.

The content of Lemma~\ref{tom} in Case~2 is then that this bundle
admits $k$ disjoint sections, one of which is the diagonal. The proof
given above was chosen because it requires no conditions on the genus
$g$ of $\Sigma$. In the special case when $k|(g-1)$, another
construction is as follows. Let $\sigma\colon \Sigma\to \Sigma$
generate a free action of $\BZ/k\BZ$ on $\Sigma$; then the graphs
$\Delta=\Gamma_{\text{id}}, \Gamma_\sigma,
\Gamma_{\sigma^2},\ldots,\Gamma_{\sigma^{k-1}}$ give $k$ disjoint
sections of $\Sigma\times\Sigma$.\bigskip

\noindent \textbf{Fiberwise branched covers.} In Case~3, we exploit
the connection between $\Map(\Sigma)$ and $\Map(S,\z)$, where
$S=\Sigma/\langle\tau\rangle$ and $\z$ is the image of the fixed
points of $\tau$. For surface bundles, this corresponds to passing to a
fiberwise branched cover, as follows; we allow the order of $\tau$ to
be any $k\geq 3$. If $S\to E\to B$ is a surface bundle with $n$
disjoint sections $\sigma_1,\ldots,\sigma_n\colon B\to E$, the union
of the sections gives a (disconnected) codimension 2 subspace of
$E$. Depending on the bundle and sections, $E$ may admit a cyclic
branched cover $\widetilde{E}\to E$ of order $k$, branched over the
sections $\sigma_i$; in this case $\widetilde{E}$ becomes a
$\Sigma$--bundle $\Sigma\to \widetilde{E}\to B$. The action of $\tau$
on $\Sigma$ then corresponds to the order--$k$ automorphism $\CT\colon
\widetilde{E}\to\widetilde{E}$ generating the deck transformations of
the branched cover $\widetilde{E}\to E$. The observation above that
$C(T)/\langle T\rangle$ has finite index in $\Map(S,\z)$ becomes here
the following fact: even if $E$ does not admit such a branched cover,
there is always some finite cover $B'\to B$ so that the pullback
bundle $S\to E'\to B'$ admits a cyclic branched cover, branched over
the preimages in $E'$ of the sections $\sigma_i$.\\

Combining this construction with the choice of sections
$\Gamma_{\sigma^i}\subset \Sigma\times \Sigma$ recovers the classical
example of Kodaira \cite{Kodaira} and Atiyah \cite{Atiyah} . Their
surface bundle is constructed as follows: start with a surface $S$
admitting a free action of $\BZ/k\BZ$ generated by $\sigma$. The
bundle $S\times S\to S$ does not admit a branched cover branched over
the union of the sections $\Gamma_{\sigma^i}$. However, taking $\pi\colon
S'\to S$ to be the cover corresponding to the kernel of $\pi_1(S)\to
H_1(S)\to H_1(S;\BZ/k\BZ)$, the pullback $S'\times S\to S'$ does admit
a branched cover $M_k\to S'\times S$ of order $k$, branched over the
union of the sections $\Gamma_{\sigma^i\circ \pi}$. Composing with the
projection $S'\times S\to S'$ gives a bundle $\Sigma\to M_k\to S'$,
where the fiber $\Sigma$ is a branched cover of the original fiber $S$
of order $k$, branched over $k$ points. (Note that the manifold $M_k$
fibers over a surface in two different ways; the fibering considered
here is that of the original authors.)

Aside from the choice of sections, these steps correspond exactly to the
considerations above, and so the results of Case~3 apply identically
to this case, giving the following theorem:
\begin{named}{Theorem \ref{AK}}
  When $k\geq 3$, the Atiyah--Kodaira bundle $\Sigma\to M_k\to S'$
  admits no flat connection invariant under the order--$k$ deck
  transformation $\CT\colon M_k\to M_k$.
\end{named}

The surface group $\pi_1(S',z)\subset \Map(\Sigma)$ singled out in the
previous section is the monodromy of this surface bundle. We remark
that by returning to the choice of sections considered in Case~3, the
same theorem is obtained for the surface bundles constructed by
Gonz{\'a}lez-D{\'{\i}}ez and Harvey in \cite{GD-H}.\\

We now sketch a description of Morita's $m$--construction; this is a
generalization of the construction of Kodaira and Atiyah, used by
Morita in \cite{Morita} to give the original proof of Morita's
theorem. Roughly, the $m$--construction begins with a surface bundle
over a manifold of dimension $n$ satisfying certain conditions,
then modifies it by pulling back along covers of the base, covers and
branched covers of the fiber, and the bundle projection itself; the
result is another surface bundle whose base has dimension $n+2$.

More precisely, given an admissible surface bundle $s\to E\to B$,
first pull back to the total space to obtain a bundle over $E$ with
fiber $s$; this bundle naturally admits a ``diagonal''
section. Possibly passing to a finite cover of the base, we may take a
fiberwise cover, obtaining a new bundle with fiber $S$, where $S\to s$
is a cover with deck transformation group $\BZ/m\BZ$. As discussed
above, combining the ``diagonal'' section with this $\BZ/m\BZ$--action
yields $m$ disjoint sections of this $S$--bundle. Again possibly
passing to a finite cover of the base, we may take a fiberwise
branched cover, yielding a bundle $\Sigma\to \widetilde{E}\to E'$,
where $\Sigma\to S$ is a cyclic branched cover of order $m$ branched
at $m$ points. Note that the deck transformation
$\CT\colon\widetilde{E}\to \widetilde{E}$ of this cyclic branched
cover has order $m$.

Fixing a single fiber of the original bundle $s\to E\to B$ and
following through this construction, we see that the preimage of this
fiber in $\widetilde{E}$ gives an Atiyah--Kodaira bundle $\Sigma\to
M_m\to S'$ inside $\Sigma\to \widetilde{E}\to E'$. Thus we have the
following consequence of Theorem~\ref{AK}.
\begin{sat}\label{m-const}
  When $m\geq 3$, given any admissible bundle $s\to E\to B$, the
  $\Sigma$--bundle $\Sigma\to \widetilde{E}\to E'$ resulting from
  Morita's $m$--construction admits no flat connection invariant under
  the order--$m$ deck transformation $\CT\colon \widetilde{E}\to
  \widetilde{E}$.
\end{sat}
For comparison, the corresponding form of Morita's theorem is as
follows.
\begin{sat}[Morita's Theorem]
  There exists a bundle $s\to E^6\to B^4$ so that the $\Sigma$--bundle
  $\Sigma\to \widetilde{E}^8\to E'^6$ resulting from Morita's
  $m$--construction admits no flat connection.
\end{sat}

\bigskip

{\small 
\noindent Mladen Bestvina, Department of Mathematics, University of Utah

\noindent \texttt{bestvina@math.utah.edu}

\bigskip

\noindent Thomas Church, Department of Mathematics, Stanford University

\noindent \texttt{church@math.stanford.edu}

\bigskip

\noindent Juan Souto, Mathematics Department, University of British Columbia

\noindent \texttt{jsouto@math.ubc.ca}}


\begin{thebibliography}{9}

\bibitem{Atiyah} M. F. Atiyah, \emph{The signature of fibre-bundles},
  in Global Analysis (Papers in Honor of K. Kodaira), 
  University of Tokyo Press (1969).

\bibitem{BLM} J. Birman, A. Lubotzky, and J. McCarthy, \emph{Abelian
    and solvable subgroups of the mapping class groups}, Duke Math. J.
    50 (1983).

\bibitem{Cantat} S. Cantat and D. Cerveau, \emph{Analytic actions of
    mapping class groups on surfaces}, J. Topol. 1 (2008). 

\bibitem{EarleEells} C. J. Earle and J. Eells, \emph{The
    diffeomorphism group of a compact Riemann surface},
  Bull. Amer. Math. Soc. 73 (1967).

\bibitem{EarleSchatz} C. J. Earle and A. Schatz, \emph{Teichm\"uller
    theory for surfaces with boundary}, J. Diff. Geom. 4
  (1970). 

\bibitem{F-H} J. Franks and M. Handel, \emph{Global fixed points for
    centralizers and Morita's Theorem}, Geom. Topol. 13 (2009). 

\bibitem{G} W. Goldman, \emph{Discontinuous groups and the Euler
  class}, Doctoral dissertation, Univ. of California, Berkeley, 1980.

\bibitem{GD-H} G. Gonz{\'a}lez-D{\'{\i}}ez and W. J. Harvey,
  \emph{Surface groups inside mapping class groups}, Topology 38
  (1999).

\bibitem{Kerckhoff} S. Kerckhoff, \emph{The Nielsen realization
    problem}, Ann. of Math. 117 (1983).

\bibitem{Kodaira} K. Kodaira, \emph{A certain type of irregular
    algebraic surfaces}, J. Analyse Math. 19 (1967).

\bibitem{KM} D. Kotschick, S. Morita, \emph{Signatures of foliated
    surface bundles and the symplectomorphism groups of surfaces},
  Topology 44 (2005).

\bibitem{Markovic} V. Markovi\'{c}, \emph{Realization of the mapping class
    group by homeomorphisms}, Invent. Math. 168 (2007).

\bibitem{Markovic-Saric} V. Markovi\'{c} and D. \v{S}ari\'{c}, \emph{The mapping
    class group cannot be realized by homeomorphisms}, preprint (2008).

\bibitem{Milnor} J. Milnor, \emph{On the existence of a connection
    with curvature zero}, Comment. Math. Helv. 32 (1958).

\bibitem{MS} J. Milnor and J. Stasheff, \emph{Characteristic classes},
  Annals of Mathematics Studies, No. 76. Princeton University Press
  (1974).

\bibitem{Morita} S. Morita, \emph{Characteristic classes of surface
    bundles}, Invent. Math. 90 (1987).

\bibitem{MoritaBook} S. Morita, \emph{Geometry of characteristic
    classes}, translated from the 1999 Japanese original, Translations
  of Mathematical Monographs, 199, Iwanami Series in Modern
  Mathematics, American Mathematical Society, Providence, RI (2001).

\bibitem{Nielsen} J. Nielsen, \emph{Abbildungsklassen endlicher Ordnung}
(German), Acta Math.  75  (1943).

\bibitem{Wood} J. Wood, \emph{Bundles with totally disconnected
    structure group}, Comment. Math. Helv. 46 (1971)

\bibitem{ZVC} H. Zieschang, E. Vogt, and H.-D. Coldewey, \emph{Fl\"{a}chen
    und ebene diskontinuierliche Gruppen} (German), Lecture Notes in
  Mathematics, Vol. 122 Springer-Verlag, Berlin-New York (1970)

\end{thebibliography}
\end{document}